\documentclass[11pt]{article}
\usepackage{amsmath,amssymb,amsthm,stmaryrd}
\title{A triple construction for Lie bialgebras}
\author{Jan E. Grabowski}

\date{21st January, 2005}

\newtheorem{theorem}{Theorem}[section]
\newtheorem*{theorem*}{Theorem}
\newtheorem{proposition}[theorem]{Proposition}
\newtheorem{lemma}[theorem]{Lemma}
\newtheorem{corollary}[theorem]{Corollary}
\newtheorem{definition}[theorem]{Definition}

\newcommand{\ad}[2]{{\mathrm{ad}}_{#1}(#2)}
\newcommand{\bad}[2]{{\underline{\mathrm{ad}}}_{#1}(#2)}
\newcommand{\bLie}[1]{\underline{\Lie{#1}}}
\let\chisave\chi
\renewcommand{\chi}{{%
 \mathchoice{\raisebox{0.25ex}{$\displaystyle\chisave$}}
            {\raisebox{0.2ex}{$\textstyle\chisave$}}
            {\raisebox{0.2ex}{$\scriptstyle\chisave$}}
            {\raisebox{0.1ex}{$\scriptscriptstyle\chisave$}}}}
\newcommand{\complex}{\ensuremath \mathbb{C}}
\newcommand{\dbos}[3]{\Lie{#1} \rtimesdot \Lie{#2} \ltimesdot \Lie{#3}}
\newcommand{\dbxcosum}{\mathrel{\blacktriangleright\joinrel\mkern-3mu\blacktriangleleft}}
\newcommand{\dbxsum}{\bowtie}
\newcommand{\dbxsumrbos}{\mathrel{\rhd\joinrel\mkern-11mu\lhd\joinrel\mkern-15mu\cdot}}
\newcommand{\defeq}{\stackrel{\mathrm{def}}{=}}
\newcommand{\dform}{\mathrm{d}}
\newcommand{\dsum}{\ensuremath{ \oplus}}
\newcommand{\dual}[1]{\ensuremath {#1}^{*}}
\newcommand{\eg}{e.g.\ }
\renewcommand{\epsilon}{\varepsilon}
\newcommand{\href}[2]{#2}
\newcommand{\id}{\ensuremath \mbox{\textup{id}}}
\newcommand{\ie}{i.e.\ }

\newcommand{\inv}[1]{\ensuremath {#1}^{-1}}
\newcommand{\ip}[2]{\ensuremath \lgen\;\!#1,#2\;\!\rgen}
\newcommand{\iipp}[2]{\mathopen{\ll}\;\!#1,#2\;\!\mathclose{\gg}}
\newcommand{\iso}{\ensuremath \cong}
\newcommand{\itund}[1]{{#1}{\!\!\!\!\!\underline{\ \,}}}
\newcommand{\laction}{\triangleright}
\newcommand{\Lbbracket}[2]{\ensuremath{ \llbracket\, #1 , #2\, \rrbracket}}
\newcommand{\Lbracket}[2]{\ensuremath{ [\, #1 , #2\, ]}}
\newcommand{\lgen}{\ensuremath \mathopen{<}}
\newcommand{\Lie}[1]{\ensuremath{\mathfrak{#1}}}
\renewcommand{\ltimes}{\mathrel{\rhd\joinrel\mkern-6mu<}}
\newcommand{\ltimesdot}{\mathrel{\cdot\joinrel\mkern-14.3mu\rhd\mkern-8.7mu<}}
\newcommand{\modcat}[1]{\ensuremath \!_{#1}\mathcal{M}}
\newcommand{\op}[1]{{#1}^{\mbox{\scriptsize \textup{op}}}}
\renewcommand{\phi}{\varphi}

\newcommand{\reals}{\ensuremath \mathbb{R}}
\newcommand{\rgen}{\ensuremath \mathclose{>}}
\renewcommand{\rtimes}{\mathrel{>\joinrel\mkern-6mu\lhd}}
\newcommand{\rtimesblack}{\mathrel{>\joinrel\mkern-1mu\blacktriangleleft}}
\newcommand{\rtimesdot}{\mathrel{>\joinrel\mkern-6mu\lhd\mkern-7.2mu\cdot}}
\newcommand{\setspan}[2]{\mathrm{span}_{#1}\{ #2 \}}
\newcommand{\smun}[2]{{#1}_{\scriptscriptstyle #2}}
\newcommand{\tensor}{\ensuremath \otimes}
\newcommand{\webtilde}{\small{$\sim\!$} \normalsize}

\input cyracc.def

\numberwithin{equation}{section}

\begin{document}

\maketitle

\begin{abstract}

\noindent We study the triple of a quasitriangular Lie bialgebra  as a natural extension of the Drinfel\cprime d double.  The triple is itself a quasitriangular Lie bialgebra.  We prove several results about the algebraic structure of the triple, analogous to known results for the double.  Among them, we prove that in the factorisable case the triple is isomorphic to a twisting of $\Lie g\oplus \Lie g\oplus \Lie g$ by a certain cocycle.  We also consider real forms of the triple and the triangular case.
\end{abstract}

\section{Introduction}

The study of Lie bialgebras (\cite{Drinfeld2}) is now well established as an infinitesimalisation of the notion of a quantum group or Hopf algebra.  A Lie bialgebra consists of a Lie algebra $\Lie{g}$ and a compatible Lie cobracket $\delta:\Lie{g} \to \Lie{g} \tensor \Lie{g}$ that controls the tensor product of representations.  Equivalently, the cobracket $\delta$ induces a bracket on the dual $\dual{\Lie{g}}$.  Moreover, Lie bialgebras exponentiate geometrically to Poisson--Lie groups with Poisson bracket linearizing to $\delta$ and have been of considerable interest to Poisson and symplectic geometers. We recall further that a Lie bialgebra is quasitriangular if $\delta$ is the appropriate coboundary of an element $r\in \Lie{g} \tensor \Lie{g}$ obeying the classical Yang--Baxter equations.  A canonical example of this type exists for all simple $\Lie{g}$, leading to the Drinfel\cprime d--Sklyanin Poisson bracket on the associated Lie group.

Probably the most important of all Lie bialgebra constructions is the Drinfel\cprime d double that associates to a Lie bialgebra a quasitriangular one $D(\Lie{g})=\Lie{g} \dbxsum \op{\dual{\Lie{g}}}$. Here it is presented as a double cross sum of $\Lie{g}$ and its dual acting on each other. When $\Lie{g}$ is quasitriangular one also knows that $D(\Lie{g})$ is isomorphic to a `bosonisation' $\Lie{g} \ltimesdot \op{\dual{\bLie{g}}}$ as well as a cocycle twist of $\Lie{g} \dsum \Lie{g}$ when $\Lie{g}$ is factorisable (when the symmetric part of $r$ is non-degenerate). Here $\bLie{g}$ denotes the braided-Lie bialgebra associated to $\Lie{g}$ (\cite{BraidedLie}), a theory which infinitesimalises the theory of braided groups, that is, Hopf algebras in braided categories. Indeed these results about the Drinfel\cprime d double infinitesimalise results about the Drinfel\cprime d quantum double in terms of braided groups and Hopf algebra twisting respectively.

There are, however, some defects of the Drinfel\cprime d double and these are solved in a canonical `triple' construction $T(\Lie{g})$ which we formulate and study here. First of all, the Drinfel\cprime d double does not respect the Cartan decomposition of Lie groups and Lie algebras into positive roots, Cartan and negative roots, so cannot be used to construct them directly. (Drinfel\cprime d here, and in the quantum version to construct quantum groups, uses the double of the Borel subalgebra and then has to identify the two Cartan subalgebras via a quotient.)  Related to this, the double is a special case of constructions (as above) which do not in general yield quasitriangular Lie bialgebras so this key property is again imposed by hand. By contrast our triple  

\[ T(\Lie{g}) \defeq\dbos{\underline{g}}{g}{\op{\dual{\underline{g}}}}\]

is a canonical example of a general triple product or `double
bosonisation' construction $\dbos{b}{g}{\dual{b}}$ (\cite{BraidedLie})
which is always quasitriangular and which respects the Cartan
decomposition of simple Lie algebras (and of quantum groups in the
quantum case). The special case in which $\Lie{b}=\bLie{g}$ and the
actions are (co)adjoint is canonical and provides our natural
extension of the double. After formulating $T(\Lie{g})$ we prove that
it is indeed an extension of the Drinfel\cprime d double and prove
several theorems about it.  Our principal result (Theorem~\ref{twist})
is the isomorphism of $T(\Lie{g})$ with $(\Lie{g} \dsum \Lie{g} \dsum
\Lie{g})_{\chi}$, the twisting of the direct sum bialgebra by a
cocycle.  This is precisely analogous to that for the double mentioned
above.  Another result is that $T(\Lie{g})$ is isomorphic to $\bLie{g}
\dbxsumrbos D(\Lie{g})$ (Corollary~\ref{dbxsumrbos}), i.e. a double
cross sum as a Lie algebra and semidirect as a Lie coalgebra.

Apart from its internal structure, some secondary motivations are as follows. First of all, as mentioned above and described in \cite{BraidedLie}, double-bosonisation may be used to construct Lie (bi-)algebras in an inductive manner: given $\Lie{g}$ we adjoin more negative and positive roots $\Lie{b}$, $\dual{\Lie{b}}$ respectively to obtain another quasitriangular Lie bialgebra.  One may use this to attempt to understand the construction and classification of simple (complex) Lie algebras. For example, one could ask why there is no $E_{9}$ arising from $\Lie{g}=E_{8}$, for which we would need to consider braided-Lie bialgebras in the category of modules of $E_{8}$.  Since the smallest non-trivial representation of $E_{8}$ is the adjoint one, one would be led naturally to the quasitriangular Lie bialgebra  $\underline{E_{8}} \rtimesdot E_{8} \ltimesdot \underline{E_{8}}=T(E_{8})$ as the smallest candidate for $E_{9}$ from this point of view. Our Theorem~\ref{twist} tells us that this is not simple as a Lie algebra, as expected, but note that it is non-trivial and not a direct sum as a Lie \emph{bi}algebra. 

Secondly, Lie bialgebras (and the Drinfel\cprime d double in particular) play a central role in mathematical physics---in the theory of integrable systems and more recently in string theory.  For example, both the Lie bialgebra $D(\Lie{g})$ and its corresponding Poisson-Lie group $D(G)$ occur when one considers Poisson-Lie $T$-duality (\cite{Klimcik}) in $\sigma$-models.  The algebraic structure behind this has been analysed in \cite{BeggsMajid}, making central use of the cotwist theorem for the Drinfel\cprime d double mentioned above.  Beggs and Majid also generalise $T$-duality to more general double cross sums.  It seems likely then that our extension $T(\Lie{g})$, which has both the cotwist property and can be represented as a double cross sum, should be the natural basis for an extension of Poisson-Lie $T$-duality, possibly in the context of higher order dualities based on a `triple product' factorisation rather than a binary factorisation. Likewise for other models in mathematical physics where the Drinfel\cprime d double is used. 

Finally, all our results should be the semi-classical version of, and provide insight into, quantum group versions of similar constructions.  The analogous general double-bosonisation is known (\cite{DoubleBos}) and the special case which we would call 
 \[ T(H) \defeq \itund{H}\rtimesdot H \ltimesdot \op{\dual{\itund{H}\;}}\]  
is again a canonical example (using quantum group (co)adjoint actions). Its particular structure has not been studied but it is quasitriangular (from the general theory) at least in the finite-dimensional case, and is an extension of the celebrated Drinfel\cprime d quantum double $D(H)$. In terms of applications,
the triple $T(H)$ can be expected to extend the role of the double $D(H)$. For example, in non-commutative differential geometry the bicovariant differential calculi on a Hopf algebra $H$ were classified by Woronowicz (\cite{Woron}) effectively in terms of the representations of $D(H)$.   In \cite{MajDiffCalc}, Majid notes that braided bicovariant differential calculi on braided groups $\itund{H}\:$ are classified in an entirely analogous way precisely by the appropriate double-bosonisation, which turns out to be $T(H)$.  If one could prove a cotwisting theorem for $T(H)$---that as an algebra it is the tensor product of three copies of $H$ (this is suggested by our Theorem~\ref{dirsum} for $T(\Lie g)$)---one would then be able to classify such braided differential calculi, for example. Such a theorem, if true, appears to be rather non-trivial to prove. 

The structure of this paper is as follows.  We begin with a section of preliminaries (Section~\ref{prelims}), recalling the definitions of Lie bialgebras and the Drinfel\cprime d double.  We also recall braided-Lie bialgebras and state the general theorem defining double-bosonisation.  Section~\ref{triple} studies $T(\Lie{g})$ and contains our first simplifications of its algebraic structure in the factorisable case.  Section~\ref{structure} contains our main result that the Lie algebra structure is isomorphic to  the direct sum $\Lie{g} \dsum \Lie{g} \dsum \Lie{g}$ (Theorem~\ref{dirsum}) and a cocycle twisting as a Lie bialgebra  (Theorem~\ref{twist}).  Note that these results do not imply that the triple $T(\Lie{g})$ is trivial any more than for $D(\Lie g)$.  Next, in Section~\ref{dbreln}, we give several theorems relating the triple to the Drinfel\cprime d double, as an extension, in various ways.

All these general results are over any field of characteristic not 2.  In Section~\ref{realforms}, we work over $\complex$ and examine the half-real forms of the triple.   A half-real form is a choice of basis of a complex bialgebra with real Lie algebra structure constants and imaginary coalgebra structure constants.  We prove that if $(\Lie{u},r)$ is a half-real form of $(\Lie{g},r)$ then $T(\Lie{u})$, defined to be the appropriate double-bosonisation, is a half-real form of $T(\Lie{g})$, under a natural reality assumption on $r$, and that $T(\Lie{u})$ has its quasitriangular structure also of this real type.

Finally, for completeness, we discuss the case when our input bialgebra is triangular---the opposite extreme from factorisable.  Whereas factorisability is the property that the symmetric part $2r_{+}$ of the quasitriangular structure defines an isomorphism of $\Lie{g}$ with $\dual{\Lie{g}}$, triangularity is the case where this symmetric part is zero.  We examine this situation in Section~\ref{triang}.

It is our opinion that the results we have obtained here and the applications described above strengthen the argument that among the variety of double cross products and coproducts, the study of double-bosonisation in particular is well justified.  Much further work is clearly necessary, in particular for the quantum triple, and we hope this paper will stimulate that.

The author would like to thank Shahn Majid for much help and guidance.  The author also gratefully acknowledges financial support from the EPSRC.

\section{Preliminaries}\label{prelims}

Throughout, unless otherwise stated, we work over a field $k$ of characteristic not 2.  We use $\tau$ to mean the tensor product flip map, \eg \[ \tau:V\tensor W \to W\tensor V,\ \tau(v\tensor w)=w\tensor v\ \mbox{for all}\ v\in V,\ w\in W, \] on any appropriate pair of vector spaces.  The adjoint action of $\Lie{g}$ on itself can be extended naturally to tensor products as follows.  For $x,y,z\in \Lie{g}$, \[ \ad{x}{y\tensor z}=\ad{x}{y}\tensor z +y\tensor \ad{x}{z}. \]  We will use this throughout without further comment.  We use the term ad-invariant in the obvious way.

We adopt the Sweedler notation for elements of tensor products.  That is, we use upper or lower parenthesized indices to indicate the placement in the tensor product, \eg $\sum a_{(1)} \tensor a_{(2)} \tensor a_{(3)} \in \Lie{g} \tensor \Lie{g} \tensor \Lie{g}$.  We will usually drop the summation sign, as in the Einstein convention.

The definition of a Lie bialgebra is originally due to Drinfel\cprime d (\cite{Drinfeld1}, \cite{Drinfeld2}).  The idea is the same as that for Hopf algebras, where we have two structures dual to each other, compatible in a natural way.  It is worth commenting that Lie bialgebras form a richer class than Lie algebras: the choice of the cobracket, the dual structure to the bracket, is not usually unique.

\begin{definition}[{\cite{Drinfeld1}}] A Lie bialgebra is $(\Lie{g},\Lbracket{\ }{\ },\delta)$ where 
\begin{enumerate} 
\item $(\Lie{g},\Lbracket{\ }{\ })$ is a Lie algebra, 
\item $(\Lie{g},\delta)$ is a Lie coalgebra, that is, $\delta:\Lie{g}\to \Lie{g}\tensor \Lie{g}$ satisfies 
	\[ \begin{array}{lll} \delta+\tau\circ \delta = 0 & & \mbox{\textup{(anticocommutativity)}} \\
			(\delta \tensor \id)\circ \delta +\mbox{\textup{cyclic}}=0 & & 
							\mbox{\textup{(co-Jacobi identity)}}
	\end{array} \]  (Here,\ \textup{``cyclic''}\ refers to cyclical rotations of the three tensor product factors in $\Lie{g}\tensor \Lie{g} \tensor \Lie{g}$.)
\item we have a cohomological compatibility condition: $\delta$ is a 1-cocycle in $Z_{\mathrm{ad}}^{1}(\Lie{g},\Lie{g}\tensor\Lie{g})$.  Explicitly, \[ \delta(\Lbracket{x}{y})=\ad{x}{\delta y}-\ad{y}{\delta x}. \] 
\end{enumerate}
\end{definition}

\noindent Examining this definition, we see that if $\Lie{g}$ is a finite-dimensional Lie bialgebra, then $(\dual{\Lie{g}},\dual{\delta},\dual{\Lbracket{\ }{\ }})$ is also a finite-dimensional Lie bialgebra.  Here, $\dual{\delta}$ and $\dual{\Lbracket{\ }{\ }}$ are the bracket and cobracket, respectively, given by dualisation.

In many of the natural cases we wish to consider, the cobracket $\delta$ arises as the coboundary of an element $r\in \Lie{g}\tensor\Lie{g}$.  Explicitly, $\delta x=\ad{x}{r}$ for all $x\in \Lie{g}$.  Imposing the further conditions that $r$ satisfies the classical Yang--Baxter equation and has ad-invariant symmetric part, we say that $(\Lie{g},r)$ is a quasitriangular Lie bialgebra.

The classical Yang--Baxter equation, in the Lie setting, is \[ \Lbbracket{r}{r} \stackrel{\mathrm{def}}{=} \Lbracket{\smun{r}{12}}{\smun{r}{13}} + \Lbracket{\smun{r}{12}}{\smun{r}{23}} + \Lbracket{\smun{r}{13}}{\smun{r}{23}} = 0. \]  Here, writing $r=r^{(1)} \tensor r^{(2)}$, $\smun{r}{12}=r^{(1)}\tensor r^{(2)} \tensor 1$, etc., with the indices showing the placement in the triple tensor product $\Lie{g} \tensor \Lie{g} \tensor \Lie{g}$.  The bracket is taken in the common factor, \eg $\Lbracket{\smun{r}{12}}{\smun{r}{13}}=\Lbracket{r^{(1)}}{{r'}^{(1)}} \tensor r^{(2)} \tensor {r'}^{(2)}$ with $r'$ a second copy of $r$.  The bracket $\Lbbracket{\ }{\ }$ is the Schouten bracket, the natural extension of the bracket to these tensor spaces.  We can take $\Lbbracket{r}{s}$ in the above definition by replacing $r'$ with $s$.

We consider the element $\Lie{q}\in \Lie{g}\tensor \Lie{g}$ where $\Lie{q}=r+\smun{r}{21}$, \ie twice the symmetric part of $r$.  We will usually write $2r_{+}$ for $\Lie{q}$.  We can distinguish two further cases.  Firstly, if $\Lie{q}=0$ we say $(\Lie{g},r)$ is triangular.  Secondly, considering $\Lie{q}$ as a map $\dual{\Lie{g}} \to \Lie{g}$, if this map is surjective, we say $(\Lie{g},r)$ is factorisable.  We will use ``factorisable'' and ``triangular'' for ``factorisable quasitriangular'' and ``triangular quasitriangular'', respectively.

Note also that any quasitriangular structure $r$ gives a second for free: it is easy to check that if $r$ satisfies the conditions above then so does $-\smun{r}{21}$.  We call $-\smun{r}{21}$ the opposite quasitriangular structure to $r$.

The Drinfe\cprime d double is one of the most studied and most used Lie bialgebra constructions and it will play an equally important role in our study of the triple.

\begin{definition}[{\cite{Drinfeld2}}] Let \Lie{g} be a finite-dimensional Lie bialgebra.  The Drinfel\cprime d double, $D(\Lie{g})$, is the quasitriangular Lie bialgebra given by \begin{enumerate} \item $\Lie{g} \dsum \dual{\Lie{g}}$ as base vector space, \item Lie bracket given by $\Lie{g}$ as a sub-Lie algebra in the first part, $\op{\dual{\Lie{g}}}$ as a sub-Lie algebra in the second part and bracket between the two given by \[ \smun{\Lbracket{x}{\phi}}{D}= \phi_{(1)}\ip{\phi_{(2)}}{x} + x_{(1)}\ip{\phi}{x_{(2)}} \] for $x\in \Lie{g}$, $\phi\in \op{\dual{\Lie{g}}}$ and \item Lie cobracket given by the direct sum cobracket, \ie $\smun{\delta}{D}x=\delta x$ and $\smun{\delta}{D}\phi = \dual{\delta} \phi$ where $\delta$ and $\dual{\delta}$ are the cobrackets on $\Lie{g}$ and $\dual{\Lie{g}}$ respectively. \end{enumerate} \end{definition}

\noindent Note that the above coalgebra structure is given by taking $r=\sum f^{a} \tensor e_{a}$ with $\{e_{a}\}$ a basis of $\Lie{g}$ and $\{f^{a}\}$ a dual basis.  We see that the double is always quasitriangular.

We have a number of different realisations of the double by other constructions.  Firstly, we have the double cross sum Lie bialgebra, $\Lie{g} \dbxsum \op{\dual{\Lie{g}}}$.  Here the two parts act on each other and the cross bracket is obtained from this: we think of a simultaneous two-way semidirect sum.  The coalgebra is a direct sum.  The actions here are the coadjoint action, namely $\phi \laction x = x_{(1)}\ip{\phi}{x_{(2)}}$.  See \cite[Section 8.3]{FQGT} for more details and references.  A second is given by considering a twisted structure, as described below.  We will recall a third, the most relevant as the inspiration for the definition of the triple, at the start of Section~\ref{triple}. 

We can obtain many non-trivial structures by twisting simple ones using cohomology.  In our setting of quasitriangular Lie bialgebras we will only focus on twisting cobrackets, or equivalently, quasitriangular structures.  To twist a quasitriangular structure $r$, we replace $r$ by $r+\chi$, where $\chi\in \Lie{g}\tensor \Lie{g}$ satisfies $\Lbbracket{r}{\chi}+\Lbbracket{\chi}{r}+\Lbbracket{\chi}{\chi}=0$ and for all $\xi\in \Lie{g}$, $\ad{\xi}{\chi + \smun{\chi}{21}}=0$.  Here, $\Lbbracket{\ }{\ }$ is the Schouten bracket as described above.  One can check that $(\Lie{g},\Lbracket{\ }{\ },r+\chi)$ is indeed again a quasitriangular Lie bialgebra.  We remark that the first of these conditions is equivalent to \[ \ad{\xi}{(\id\tensor \delta)\chi + \mathrm{cyclic} + \Lbbracket{\chi}{\chi}}=0 \] for all $\xi \in \Lie{g}$.  So, in particular, we can satisfy the twisting requirements by choosing $\chi$ such that $\chi$ is symmetric and has $(\id\tensor \delta)\chi + \mathrm{cyclic} + \Lbbracket{\chi}{\chi}=0$.

Then we can describe the Drinfel\cprime d double as a twisting, as follows.

\begin{theorem}[{\cite{STS}}]\label{dbtwist} Let $(\Lie{g}, \Lbracket{\ }{\ },r)$ be a quasitriangular Lie bialgebra.  There is a quasitriangular Lie bialgebra $\Lie{g} \dbxcosum \Lie{g}$ given by twisting $\Lie{g} \dsum \Lie{g}$ by \[ \chi = \smun{r}{LR}-\tau(\smun{r}{LR})=(r^{(1)} \dsum 0) \tensor (0 \dsum r^{(2)}) - (0 \dsum r^{(2)}) \tensor (r^{(1)} \dsum 0). \]  Moreover, there is a homomorphism $D(\Lie{g}) \to \Lie{g} \dbxcosum \Lie{g}$ of Lie bialgebras which is an isomorphism when $\Lie{g}$ is factorisable. \end{theorem}

\noindent The quasitriangular structure on $\Lie{g} \dsum \Lie{g}$ we start with in this theorem is $r \dsum -\smun{r}{21}$, \ie the direct sum, taking the opposite quasitriangular structure on the second factor.

We now consider the braided version of Lie bialgebras, as defined by Majid in \cite{BraidedLie}.  Here we consider the module category $\modcat{\Lie{g}}$ of a quasitriangular Lie bialgebra $\Lie{g}$ and objects in this category possessing a $\Lie{g}$-covariant Lie algebra structure.  Following the line suggested by the theory of braided groups, we associate to these objects a braiding-type map generalising the usual flip.  If $\Lie{b}$ is a $\Lie{g}$-covariant Lie algebra in the category $\modcat{\Lie{g}}$, we define the infinitesimal braiding of $\Lie{b}$ to be the operator $\psi:\Lie{b}\tensor \Lie{b} \to \Lie{b}\tensor\Lie{b}$, $\psi(a\tensor b)=2r_{+} \laction (a\tensor b-b\tensor a)$ where $\laction$ is the left action of $\Lie{g}$ on $\Lie{b}$.  In fact, $\psi$ is a 2-cocycle in $Z_{\mathrm{ad}}^{2}(\Lie{b},\Lie{b}\tensor \Lie{b})$.

\pagebreak
\begin{definition}[{\cite{BraidedLie}}]\label{blba} A braided-Lie bialgebra $(\Lie{b},\Lbracket{\ }{\ }_{\Lie{b}},\underline{\delta})$ is an object in $\modcat{\Lie{g}}$ satisfying the following conditions.  \begin{enumerate} \item $(\Lie{b},\Lbracket{\ }{\ }_{\Lie{b}})$ is a $\Lie{g}$-covariant Lie algebra in the category. \item $(\Lie{b},\underline{\delta})$ is a $\Lie{g}$-covariant Lie coalgebra in the category. \item $\dform \underline{\delta}=\psi$. \end{enumerate} \end{definition}

\noindent We can now state the theorem which provides the construction we use in this paper.  Let $\Lie{g}$ be a quasitriangular Lie bialgebra.

\begin{theorem}[{\cite[Theorem 3.10]{BraidedLie}}]\label{dbos} For dually paired braided-Lie bialgebras $\Lie{b},\Lie{c}\in \modcat{\Lie{g}}$ the vector space $\Lie{b} \dsum \Lie{g} \dsum \Lie{c}$ has a unique Lie bialgebra structure $\dbos{b}{g}{\op{c}}$, the double-bosonisation, such that $\Lie{g}$ is a sub-Lie bialgebra, $\Lie{b},\Lie{\op{c}}$ are Lie subalgebras, and \[ \begin{array}{c} \Lbracket{\xi}{x}=\xi \laction x, \quad \Lbracket{\xi}{\phi}=\xi \laction \phi \\ \\  \Lbracket{x}{\phi} = x_{\underline{(1)}}\ip{\phi}{x_{\underline{(2)}}}+\phi_{\underline{(1)}}\ip{\phi_{\underline{(2)}}}{x}+2r_{+}^{(1)}\ip{\phi}{r_{+}^{(2)}\laction x} \\ \\ \delta x=\underline{\delta}x+r^{(2)}\tensor r^{(1)}\laction x-r^{(1)}\laction x\tensor r^{(2)} \\ \\ \delta \phi=\underline{\delta}\phi + r^{(2)}\laction \phi \tensor r^{(1)} - r^{(1)}\tensor r^{(2)}\laction \phi \end{array} \] \noindent for all $x\in \Lie{b}$, $\xi\in \Lie{g}$ and $\phi \in \Lie{c}$.  Here $\underline{\delta}x=x_{\underline{(1)}}\tensor x_{\underline{(2)}}$. \end{theorem}

\noindent Moreover, the double-bosonisation is always quasitriangular.  This is established in the case we will need by the following proposition.

\begin{proposition}[{\cite[Proposition 3.11]{BraidedLie}}]\label{doubleqtstr} Let $\Lie{b}\in \modcat{\Lie{g}}$ be a finite-dimensional braided-Lie bialgebra with dual $\dual{\Lie{b}}\!$.  Then the double-bosonisation $\dbos{b}{g}{\op{\dual{b}}}$ is quasitriangular with \[ r^{\mathrm{new}}=r+\sum_{a} f^a \tensor e_a \] where $\{ e_a \}$ is a basis of $\Lie{b}$ and $\{ f^a \}$ is a dual basis, and $r$ is the quasitriangular structure of $\Lie{g}$.  If $\Lie{g}$ is factorisable then so is the double-bosonisation.
\end{proposition}

We will occasionally refer to the (single) bosonisation, in particular in the following section.  It occurs as a sub-Lie bialgebra of the double bosonisation.  With the notation of Theorem~\ref{dbos} above, $\Lie{b} \rtimesdot \Lie{g}$ is the sub-Lie bialgebra of $\dbos{b}{g}{\op{c}}$ with base vector space $\Lie{b}\dsum \Lie{g}$ and Lie algebra and coalgebra structures as in that theorem.

\section{The triple of a Lie bialgebra}\label{triple}

We wish to discuss a special case of the double-bosonisation theorem, as recalled above.  The braided-Lie bialgebra that we will use is the braided structure $\bLie{g}$ naturally associated to any quasitriangular Lie bialgebra $\Lie{g}$, coming from the adjoint representation (\cite{BraidedLie}). 

\begin{definition}\label{transmut} Take $\bLie{g}$ to be the adjoint representation of $\Lie{g}$.  For the Lie bracket of $\bLie{g}$, we take the Lie bracket of $\Lie{g}$.  Clearly, this is covariant.  For the braided cobracket we take \[ \underline{\delta}x=2r_{+}^{(1)} \tensor \Lbracket{x}{r_{+}^{(2)}} \] for all $x\in \bLie{g}$.  Here, $r$ is considered as an element of $\bLie{g}\tensor \bLie{g}$ in the obvious way.  We call $\bLie{g}$ the transmutation of $\Lie{g}$. \end{definition}

\noindent For a finite-dimensional factorisable Lie bialgebra $\Lie{g}$, this braided structure is very natural.  Using the isomorphism $2r_{+}:\dual{\Lie{g}}\to \Lie{g}$ provided by the factorisability assumption, we have that $\underline{\delta}$ is equivalent to the Kirillov--Kostant Lie cobracket, which is the cobracket given precisely by dualising the Lie bracket of $\Lie{g}$.  Then we see that $\bLie{g}$ is self-dual. 

So we have the third description of the Drinfel\cprime d double, now as a single bosonisation, as follows.  We let $\Lie{g}$ be a finite-dimensional quasitriangular Lie bialgebra and $\dual{\bLie{g}}$ the dual of its transmutation.  The bosonisation $\dual{\bLie{g}} \rtimesdot \Lie{g}$ is isomorphic as a Lie bialgebra to the Drinfel\cprime d double $D(\Lie{g})$ (\cite[Example 3.9]{BraidedLie}).

In the factorisable case, the dual $\dual{\bLie{g}}$ here can of course be replaced by $\bLie{g}$.  As a Lie algebra, we have a semidirect sum by definition of the bosonisation and furthermore we can easily see that a semidirect sum can be re-diagonalised to a direct sum.  The coalgebra structure on this direct sum induced by these isomorphism is in fact precisely the one giving the double as a twisting (Theorem~\ref{dbtwist}).  Note that we also have $D(\Lie{g}) \iso \Lie{g} \ltimesdot \op{\dual{\bLie{g}}}$.  This isomorphism will be described and used in Section~\ref{dbreln}.

We now define the triple, as a double-bosonisation using the transmutation of $\Lie{g}$.  Let $\Lie{g}$ be a finite-dimensional Lie bialgebra over a field $k$ of characteristic not $2$.  

\begin{definition}\label{tripledefn} In Theorem~\ref{dbos}, set $\Lie{b}=\bLie{g}$ the transmutation of $\Lie{g}$ and $\Lie{c}=\dual{\bLie{g}}$.  We have $\bLie{g}\in \modcat{\Lie{g}}$ by $\mathrm{ad}$.  Define $T(\Lie{g})=\dbos{\underline{g}}{g}{\op{\dual{\underline{g}}}}$, as a Lie bialgebra.  \end{definition}

\noindent The ``$T$'' stands for ``triple'': we will show later the comparison with the Drinfel\cprime d double $D(\Lie{g})$.  We have the Lie bialgebra structure of $T(\Lie{g})$ given explicitly in terms of the brackets and cobrackets on $\Lie{g}$, $\bLie{g}$ and the module structures for $\bLie{g}$.  We will show that these formul\ae\ simplify.  

We now restrict to the case of $\Lie{g}$ factorisable, so that we can replace $\Lie{c}=\dual{\bLie{g}}$ in the above definition by $\Lie{c}=\bLie{g}$.  The pairing we use, $\iipp{\,}{\,}$, is the Killing form $K:\Lie{g}\tensor \Lie{g} \to k$ which pairs $\Lie{b}=\bLie{g}$ and $\Lie{c}=\bLie{g}$ as braided-Lie bialgebras.

For clarity, we will refer to the three pieces from left to right in the definition of $T(\Lie{g})=\dbos{\underline{g}}{g}{\op{\underline{g}}}$ as $\Lie{b}$, $\Lie{g}$ and $\op{\Lie{c}}$, respectively.  We will indicate the bracket and cobracket in $T(\Lie{g})$ by a subscript ``$T$'', to distinguish it from the brackets and cobrackets of the individual pieces, which will carry a subscript $\Lie{b},\Lie{g},\Lie{c}$ or $\op{\Lie{c}}$ as appropriate.

\begin{lemma}  As Lie algebras, let $\Lie{b}=\bLie{g}=\Lie{c}$.  The Lie bialgebra structure of $T(\Lie{g})=\dbos{b}{g}{\op{c}}$ is given by: \[ \begin{array}{c} \begin{array}{ll} \smun{\Lbracket{b_1}{b_2}}{T} = \Lbracket{b_1}{b_2}_{\Lie{b}} & \smun{\Lbracket{g}{b}}{T} = \alpha(\Lbracket{g}{\inv{\alpha}(b)}_{\Lie{g}}) \\ & \\ \smun{\Lbracket{g_1}{g_2}}{T} = \Lbracket{g_1}{g_2}_{\Lie{g}} & \smun{\Lbracket{g}{c}}{T} = \epsilon(\Lbracket{g}{\inv{\epsilon}(c)}_{\Lie{g}}) \\ & \\ \smun{\Lbracket{c_1}{c_2}}{T} = -\Lbracket{c_1}{c_2}_{\Lie{c}} &  \\ & \end{array} \\
\begin{array}{l} \smun{\Lbracket{b}{c}}{T} = -\Lbracket{b}{\alpha \circ \inv{\beta}(c)}_{\Lie{b}}+\Lbracket{\inv{\alpha}(b)}{\inv{\beta}(c)}_{\Lie{g}}+\Lbracket{\beta \circ \inv{\alpha}(b)}{c}_{\Lie{c}} \\ \\ \smun{\delta}{T} b = \delta b + (r^{(1)}\laction b)\tensor (\alpha(r^{(2)})-r^{(2)}) -(\alpha(r^{(2)})-r^{(2)})\tensor (r^{(1)} \laction b) \\ \\ \smun{\delta}{T} g = \delta g \\ \\ \smun{\delta}{T} c = \delta c  + (\beta(r^{(1)})-r^{(1)}) \tensor (r^{(2)} \laction c) - (r^{(2)} \laction c) \tensor (\beta(r^{(1)})-r^{(1)}) \end{array} \end{array} \] for $b,b_{i}\in \Lie{b}$, $g,g_{i}\in \Lie{g}$ and $c,c_{i}\in \Lie{c}$.  Here, $\alpha,\ \beta,\ \epsilon$ are the identity map between the pieces as detailed below.
\end{lemma}

\begin{proof}
We will use several isomorphisms between the three pieces in what follows.  We set \[ \begin{array}{ll} \alpha : \Lie{g} \to \Lie{b}=\bLie{g}, & \alpha = \id \\ \beta : \Lie{g} \to \Lie{c}=\bLie{g}, & \beta = \id \\ \gamma : \Lie{c} \to \op{\Lie{c}}, & \gamma= \id \\ \bar{\gamma} : \Lie{c} \to \op{\Lie{c}}, & \bar{\gamma} : c \mapsto -c \\ \epsilon : \Lie{g} \to \op{\Lie{c}}, & \epsilon = \gamma \circ \beta = \id \quad \mbox{and} \\ \bar{\epsilon} : \Lie{g} \to \op{\Lie{c}}, & \bar{\epsilon} = \bar{\gamma} \circ \beta = -\epsilon = -\id. \end{array} \] All of these except $\gamma$ and $\epsilon$ are Lie algebra isomorphisms.  We can now write the pairing of $\Lie{b}$ and $\Lie{c}$ explicitly as \[ \iipp{\,}{\,}:\Lie{b}\tensor \Lie{c}\to k, \quad \iipp{b}{c}=K(\inv{\alpha}(b),\inv{\beta}(c)) \] for $b\in \Lie{b}$, $c\in \Lie{c}$.
 
Firstly, $\Lie{g}$ has (unbraided) Lie bialgebra bracket and cobracket structures by assumption: we will denote these by plain brackets, $\Lbracket{\,\,}{\,}$, and $\delta$ respectively.  The braided-Lie bialgebra structure for $\bLie{g}$ is that described in Definition~\ref{transmut}.  The braided cobracket is \[ \underline{\delta}b = \alpha(2r_{+}^{(1)})\tensor \Lbracket{b}{\alpha(r_{+}^{(2)})}_{\Lie{b}}, \] for $b\in \Lie{b}=\bLie{g}$, where $r$ is the quasitriangular structure on $\Lie{g}$ and $2r_{+}$ is its symmetric part.  Finally, the braided-Lie bialgebra structure of $\Lie{c}=\bLie{g}$ is the same as that of $\Lie{b}$, so \begin{align*} \Lbracket{\,\,}{\,}_{\Lie{c}} & = \Lbracket{\,\,}{\,}_{\Lie{b}} \\ & = \Lbracket{\,\,}{\,}_{\Lie{g}}\quad \mbox{and} \\ \underline{\delta}c & = \beta(2r_{+}^{(1)})\tensor \Lbracket{c}{\beta(r_{+}^{(2)})}_{\Lie{c}} \end{align*}  

Each piece appears as a Lie subalgebra, so we need now to clarify the brackets between the pieces.  We have $\bLie{g}\in \modcat{\Lie{g}}$ by $ g \laction b = \ad{g}{b}$ for $g \in \Lie{g}$, $b\in \bLie{g}$, that is, the adjoint action.  To be even more explicit, for $g \in \Lie{g}$, $b\in \Lie{b}$ and $c \in \op{\Lie{c}}$ we have \begin{align*} g \laction b & = \Lbracket{\alpha(g)}{b}_{\Lie{b}} \\ & = \alpha\left( \Lbracket{g}{\inv{\alpha}(b)}_{\Lie{g}}\right) \\ g \laction c & = \Lbracket{\bar{\epsilon}(g)}{c}_{\Lie{c}^{\mbox{\tiny \textup{op}}}} \\ & = \bar{\epsilon}\left( \Lbracket{g}{\inv{\bar{\epsilon}}(c)}_{\Lie{g}}\right) \\ & = -\epsilon\left(\Lbracket{g}{-\inv{\epsilon}(c)}_{\Lie{g}}\right) \\ & = \epsilon\left(\Lbracket{g}{\inv{\epsilon}(c)}_{\Lie{g}}\right). \end{align*}  We set $\smun{\Lbracket{g}{b}}{T}=g \laction b$ and $\smun{\Lbracket{g}{c}}{T} = g \laction c$, so the action of $\Lie{g}$ on $\Lie{b}$ and $\op{\Lie{c}}$ is the adjoint action, with the bracket taken in $\Lie{g}$ after the appropriate isomorphism has been applied.

These brackets come from the Lie algebra structure of the (single) bosonisations $\Lie{b} \rtimesdot \Lie{g}$ and $\Lie{g} \ltimesdot \op{\Lie{c}}$ as in \cite[Theorem 3.5]{BraidedLie}.  Note that there are two minus signs which cancel, one from the reversed action in $\Lie{g} \ltimesdot \op{\Lie{c}}$ and one from the ``op'', so we just see the adjoint action of $\Lie{g}$ on $\Lie{c}$.

For the double-bosonisation, the remaining bracket is the one between $\Lie{b}$ and $\Lie{c}$.  This is given by \[ \Lbracket{b}{c} = b_{\underline{(1)}}\ip{c}{b_{\underline{(2)}}}+c_{\underline{(1)}}\ip{c_{\underline{(2)}}}{b}+2r_{+}^{(1)}\ip{c}{r_{+}^{(2)}\laction b} \] \noindent for $b\in \Lie{b}$, $c \in \Lie{c}$. Using the above definitions of $\underline{\delta}$ and the pairing $\iipp{\,}{\,}$, we obtain the following: \begin{align*} 
\smun{\Lbracket{b}{c}}{T} = & \ b_{\underline{(1)}}\iipp{c}{b_{\underline{(2)}}}+c_{\underline{(1)}}\iipp{c_{\underline{(2)}}}{b}+2r_{+}^{(1)}\iipp{c}{r_{+}^{(2)}\laction b} \\
  = & \ \alpha(2r_{+}^{(1)})\iipp{\Lbracket{b}{\alpha(r_{+}^{(2)})}_{\Lie{b}}}{c} + \beta(2r_{+}^{(1)})\iipp{b}{\Lbracket{c}{\beta(r_{+}^{(2)})}_{\Lie{c}}} \\ & \quad {} + 2r_{+}^{(1)}\iipp{\alpha(\Lbracket{r_{+}^{(2)}}{\inv{\alpha}(b)}_{\Lie{g}})}{c} \intertext{so by the ad-invariance of $2r_{+}$,} \smun{\Lbracket{b}{c}}{T} = & \ \Lbracket{\alpha(2r_{+}^{(1)})}{b}_{\Lie{b}}\iipp{\alpha(r_{+}^{(2)})}{c} + \Lbracket{\beta(2r_{+}^{(1)})}{c}_{\Lie{c}}\iipp{b}{\beta(r_{+}^{(2)})} \\ & \quad {} - \Lbracket{2r_{+}^{(1)}}{\inv{\alpha}(b)}_{\Lie{g}}\iipp{\alpha(r_{+}^{(2)})}{c} \\ = & \ \Lbracket{\alpha(2r_{+}^{(1)})K(\inv{\alpha}(\alpha(r_{+}^{(2)})),\inv{\beta}(c))}{b}_{\Lie{b}} \\ & \quad {} + \Lbracket{\beta(2r_{+}^{(1)})K(\inv{\alpha}(b),\inv{\beta}(\beta(r_{+}^{(2)}))}{c}_{\Lie{c}} \\ & \quad {} - \Lbracket{2r_{+}^{(1)}K(\inv{\alpha}(\alpha(r_{+}^{(2)})),\inv{\beta}(c))}{\inv{\alpha}(b)}_{\Lie{g}}. \intertext{Then, using the quasitriangular expression for the inverse of the Killing form, that is, $\inv{K}(\phi)=2r_{+}^{(1)}\phi(r_{+}^{(2)})$,} \smun{\Lbracket{b}{c}}{T} = & \ \Lbracket{\inv{K}(K(\inv{\beta}(c)))}{b}_{\Lie{b}} + \Lbracket{\inv{K}(K(\inv{\alpha}(b)))}{c}_{\Lie{c}} \\ & \quad {} - \Lbracket{\inv{K}(K(\inv{\beta}(c)))}{\inv{\alpha}(b)}_{\Lie{g}} \\ = & \ -\Lbracket{b}{\alpha \circ \inv{\beta}(c)}_{\Lie{b}}+\Lbracket{\beta \circ \inv{\alpha}(b)}{c}_{\Lie{c}}+\Lbracket{\inv{\alpha}(b)}{\inv{\beta}(c)}_{\Lie{g}} \\ = & \ -\Lbracket{b}{\alpha \circ \inv{\beta}(c)}_{\Lie{b}}+\Lbracket{\inv{\alpha}(b)}{\inv{\beta}(c)}_{\Lie{g}}-\Lbracket{\beta \circ \inv{\alpha}(b)}{c}_{\Lie{c^{\mbox{\tiny \textup{op}}}}} \\ = & \ -\smun{\Lbracket{b}{\alpha \circ \inv{\beta}(c)}}{T}+\smun{\Lbracket{\inv{\alpha}(b)}{\inv{\beta}(c)}}{T}-\smun{\Lbracket{\beta \circ \inv{\alpha}(b)}{c}}{T}.  \end{align*} \noindent  That is, the bracket of an element $b\in \Lie{b}$ with an element $c\in \Lie{c}$ is given by mapping $b$ and/or $c$ into each piece in turn and taking the bracket there.  If the bracket is non-zero, it has a non-zero component in each piece.

We now consider the Lie coalgebra structure.  We have that $\Lie{g}$ is a sub-Lie bialgebra, so the cobracket on an element of $\Lie{g}$ is simply $\delta$.  For an element $b\in \Lie{b}=\bLie{g}$, the unbraided cobracket structure of the double-bosonisation is \begin{align*} \smun{\delta}{T}b = & \ \underline{\delta}b+r^{(2)}\tensor (r^{(1)} \laction b) - (r^{(1)}\laction b)\tensor r^{(2)} \\ = & \ \delta b + (r^{(1)}\laction b)\tensor \alpha(r^{(2)}) -\alpha(r^{(2)})\tensor (r^{(1)} \laction b) \\ & \ + r^{(2)}\tensor (r^{(1)} \laction b) - (r^{(1)}\laction b)\tensor r^{(2)} \\ = & \ \delta b + (r^{(1)}\laction b)\tensor (\alpha(r^{(2)})-r^{(2)}) -(\alpha(r^{(2)})-r^{(2)})\tensor (r^{(1)} \laction b). \end{align*}

Since the bosonisation $\Lie{g} \ltimesdot \op{\Lie{c}}$ is taken to be that of $\op{\Lie{c}}$ in the category of $\Lie{g}$-modules with opposite infinitesimal braiding (see the proof of \cite[Theorem 3.10]{BraidedLie} for details), we have \begin{align*} \smun{\delta}{T} c = & \ \underline{\delta} c + (r^{(2)} \laction c) \tensor r^{(1)} - r^{(1)} \tensor (r^{(2)} \laction c) \\ = & \ \delta c - (r^{(2)} \laction c) \tensor \beta(r^{(1)}) + \beta(r^{(1)}) \tensor (r^{(2)} \laction c) \\ & \ + (r^{(2)} \laction c) \tensor r^{(1)} - r^{(1)} \tensor (r^{(2)} \laction c) \\ = & \ \delta c  + (\beta(r^{(1)})-r^{(1)}) \tensor (r^{(2)} \laction c) - (r^{(2)} \laction c) \tensor (\beta(r^{(1)})-r^{(1)}). \end{align*}
\end{proof}

In what follows we will want to compute the bracket on general elements of $T(\Lie{g})$ so we give this explicitly.  Our notation for general elements will be as elements of the direct sum vector space, usually writing $b$, $g$ and $c$ for elements of $\Lie{b}$, $\Lie{g}$ and $\Lie{c}$ respectively.  We will also now suppress the isomorphisms $\alpha,\ \beta,\ \ldots$ when taking brackets.

\begin{theorem}\label{bracket} Let $\Lie{g}$ be a factorisable Lie bialgebra and let $b_1 \dsum g_1 \dsum c_1$, $b_2 \dsum g_2 \dsum c_2 \in T(\Lie{g})$.  Then \begin{multline*} \smun{\Lbracket{b_1 \dsum g_1 \dsum c_1}{b_2 \dsum g_2 \dsum c_2}}{T} \\ \begin{aligned} = ( \Lbracket{b_1}{b_2}  + & \Lbracket{b_1}{g_2} - \Lbracket{b_1}{c_2} + \Lbracket{g_1}{b_2} - \Lbracket{c_1}{b_2} ) \\ {}\dsum ( & \Lbracket{b_1}{c_2} +  \Lbracket{g_1}{g_2} +  \Lbracket{c_1}{b_2}) \\ {}\dsum ( \Lbracket{c_1}{b_2} + & \Lbracket{c_1}{g_2} - \Lbracket{c_1}{c_2} + \Lbracket{g_1}{c_2} + \Lbracket{b_1}{c_2} ) \end{aligned} \end{multline*} \end{theorem}

\begin{proof} Immediate from the preceding lemma. \end{proof}

\section{The structure of $T(\Lie{g})$, $\Lie{g}$ factorisable}\label{structure}

We now investigate the structure of $T(\Lie{g})$, for $\Lie{g}$ a factorisable Lie bialgebra.  Our main results are Theorem~\ref{dirsum} and Theorem~\ref{twist}.  We see that these results are direct analogues of those for the Drinfel\cprime d double, as recalled in Section~\ref{prelims}.

\subsection{The Lie algebra structure}\label{LAtriple}

We start by examining the Lie ideals of $T(\Lie{g})$.

\begin{lemma}\label{subalg} The subspaces \begin{alignat*}{2} & I_{-} && = \setspan{k}{x \dsum (-x) \dsum 0 \mid x\in \Lie{g}} \\ & I_{0} && = \setspan{k}{x\dsum (-x) \dsum (-x) \mid x\in \Lie{g}} \\ & I_{+} && = \setspan{k}{0\dsum x\dsum x \mid x\in \Lie{g}} \end{alignat*} are Lie subalgebras of $T(\Lie{g})$. \end{lemma}

\begin{proof} we use the bracket on $T(\Lie{g})$ as given in Theorem~\ref{bracket}.  

For $I_{-}$, let $x\dsum (-x)\dsum 0$, $y\dsum (-y)\dsum 0\in I_{-}$.  Then \begin{align*} \smun{\Lbracket{x\dsum (-x)\dsum 0}{y\dsum (-y)\dsum 0}}{T} & = \left( \Lbracket{x}{y}-\Lbracket{x}{y}-\Lbracket{x}{y} \right) \dsum \left( \Lbracket{x}{y} \right) \dsum 0 \\ & = -\Lbracket{x}{y} \dsum \Lbracket{x}{y} \dsum 0 \\ & \in I_{-} \end{align*}  Similarly, for $I_{0}$, let $x\dsum (-x)\dsum (-x)$, $y\dsum (-y)\dsum (-y)\in I_{0}$.  Then \begin{align*} \smun{\Lbracket{x\dsum (-x)\dsum (-x)}{y\dsum (-y)\dsum (-y)}}{T} & = \Lbracket{x}{y} \dsum -\Lbracket{x}{y} \dsum -\Lbracket{x}{y} \\ & \in I_{0} \end{align*}  Finally, for $I_{+}$, let $0\dsum x\dsum x$, $0\dsum y\dsum y\in I_{+}$.  Then \begin{align*} \smun{\Lbracket{0\dsum x\dsum x}{0\dsum y\dsum y}}{T} & = 0 \dsum \Lbracket{x}{y} \dsum \Lbracket{x}{y} \\ & \in I_{+} \end{align*}
\end{proof}

\begin{lemma} The subalgebras $I_{-}$, $I_{0}$ and $I_{+}$ are ideals of $T(\Lie{g})$. \end{lemma}

\begin{proof} let $b\dsum g\dsum c\in T(\Lie{g})$.

 $I_{-}$: let $x\dsum (-x) \dsum 0\in I_{-}$.  Then \begin{align*} \smun{\Lbracket{x\dsum (-x)\dsum 0}{b\dsum g\dsum c}}{T} & = \left( \Lbracket{x}{b}+\Lbracket{x}{g}-\Lbracket{x}{c}-\Lbracket{x}{b} \right) \\ & \quad \quad \dsum \left( \Lbracket{x}{c}-\Lbracket{x}{g} \right) \dsum 0 \\ & = \Lbracket{x}{g-c} \dsum -\Lbracket{x}{g-c} \dsum 0 \\ & \in I_{-} \end{align*}

 $I_{0}$: let $x\dsum (-x) \dsum (-x)\in I_{0}$.  Then \begin{align*} \smun{\Lbracket{x\dsum (-x)\dsum (-x)}{b\dsum g\dsum c}}{T} & = \Lbracket{x}{b+g-c} \dsum -\Lbracket{x}{b+g-c} \\ & \quad \quad \dsum -\Lbracket{x}{b+g-c} \\ & \in I_{0} \end{align*}

 $I_{+}$: let $0\dsum x \dsum x\in I_{+}$.  Then \begin{align*} \smun{\Lbracket{0\dsum x\dsum x}{b\dsum g\dsum c}}{T} & = 0 \dsum \Lbracket{x}{b+g} \dsum \Lbracket{x}{b+g} \\ & \in I_{+} \end{align*} \end{proof}

\begin{theorem}\label{dirsum} $T(\Lie{g})$ is the direct sum of the ideals $I_{-}$, $I_{0}$ and $I_{+}$.  Hence $T(\Lie{g})$ is isomorphic to $\Lie{g}\dsum \Lie{g} \dsum \Lie{g}$ as a Lie algebra. \end{theorem}

\begin{proof} we must show that the brackets between any two of $I_{-}$, $I_{0}$ and $I_{+}$ are zero.

 $\Lbracket{I_{-}}{I_{0}}$: let $x\dsum (-x) \dsum 0 \in I_{-}$, $y\dsum (-y) \dsum (-y) \in I_{0}$.  Then  \begin{align*} \smun{\Lbracket{x\dsum (-x)\dsum 0}{y\dsum (-y)\dsum (-y)}}{T} & = \left( \Lbracket{x}{y}-\Lbracket{x}{y}+\Lbracket{x}{y}-\Lbracket{x}{y} \right) \\ & \quad \dsum \left( -\Lbracket{x}{y}+\Lbracket{x}{y} \right) \dsum \left( \Lbracket{x}{y}-\Lbracket{x}{y} \right) \\ & = 0 \dsum 0 \dsum 0. \end{align*}  Similarly, $\Lbracket{I_{+}}{I_{0}}$ and $\Lbracket{I_{+}}{I_{-}}$ are zero.  Hence $T(\Lie{g})$ is the direct sum of $I_{-}$, $I_{0}$ and $I_{+}$.

It is clear that $I_{-}$, $I_{0}$ and $I_{+}$ are each isomorphic to $\Lie{g}$ so we have the Lie algebra isomorphism $T(\Lie{g})\iso \Lie{g}\dsum \Lie{g} \dsum \Lie{g}$.  Notice, however, that the bracket on $I_{-}$ is the opposite one (Lemma~\ref{subalg}), so we can write \[ T(\Lie{g})\iso I_{+}\dsum I_{0} \dsum I_{-}\iso \Lie{g} \dsum \Lie{g} \dsum \op{\Lie{g}},\] which we recognise as the three `input' Lie algebras of $T(\Lie{g})$ with the bracket now diagonalised.  Alternatively, we can write $T(\Lie{g})\iso \Lie{g} \dsum \Lie{g} \dsum \Lie{g}$ using the usual isomorphism of a Lie algebra with its opposite.  

Explicitly, \[ \begin{array}{lllllllll} b \dsum g \dsum c & {}= & ( &  \quad \ \, \:\!\:\!\:\!\:\!  0 & {}\dsum & {}\quad \,\:\! b+g & {}\dsum & {}\quad\,\:\! b+g & ) \\ & {}\quad \quad + & ( & {}b+g-c & {}\dsum & {}-b-g+c & {}\dsum & {}-b-g+c & ) \\ & {}\quad \quad + & ( & {}\ \,\;\!\! -g+c & {}\dsum & {}\quad \quad\ \ \,\;\!\! g-c & {}\dsum & \quad \quad \ \ \, 0 & ) \end{array} \] is the decomposition of a general element of $T(\Lie{g})$ into a sum of elements in the ideals $I_{+}$, $I_{0}$ and $I_{-}$ respectively.  Then we have the two isomorphisms mentioned above. 
\begin{description} 
\item{$T(\Lie{g})\iso \Lie{g} \dsum \Lie{g} \dsum \op{\Lie{g}}$:}

 define $\theta_{1}:T(\Lie{g}) \to \Lie{g} \dsum \Lie{g} \dsum \op{\Lie{g}}$ by \[ b \dsum g \dsum c \mapsto (b+g) \dsum (b+g-c) \dsum (-g+c) \]
\item{$T(\Lie{g})\iso \Lie{g} \dsum \Lie{g} \dsum \Lie{g}$:}

 define $\theta_{2}:T(\Lie{g}) \to \Lie{g} \dsum \Lie{g} \dsum \Lie{g}$ by \[ b \dsum g \dsum c \mapsto (b+g) \dsum (b+g-c) \dsum (g-c) \]
\end{description}  It is easily checked that these are Lie algebra isomorphisms. \end{proof}

\noindent We have an immediate corollary.

\begin{corollary} The rank of $T(\Lie{g})$ is three times that of $\Lie{g}$. \qed \end{corollary}

\noindent This may also be proved independently of the above theorem by examining the possible Abelian subalgebras of $T(\Lie{g})$.  Indeed, we can use the above isomorphism to see that the Cartan subalgebra of $T(\Lie{g})$ is the direct sum of the three incarnations of the Cartan subalgebra of $\Lie{g}$.

\subsection{The Lie coalgebra structure}\label{LCAtriple}

Recall that the Lie coalgebra structure of a Lie bialgebra is completely determined by the Lie algebra and the quasitriangular structure, $r$.  From the previous section we have a Lie algebra isomorphism of $T(\Lie{g})$ with $\Lie{g} \dsum \Lie{g} \dsum \Lie{g}$.  Furthermore, double-bosonisation comes with an explicit expression for its quasitriangular structure.  This is given in Proposition~\ref{doubleqtstr}.  We will identify the image of the quasitriangular structure under the Lie algebra isomorphism and so express $T(\Lie{g})$ as a twisting by a cocycle of the direct sum structure.

\begin{theorem}\label{twist} Let $\Lie{g}$ be a factorisable Lie bialgebra.  As a Lie bialgebra, $T(\Lie{g})$ is isomorphic to $(\Lie{g} \dsum \Lie{g} \dsum \Lie{g})_{\chi}$, the twisting by \[ \chi= \smun{r}{AB}-\tau(\smun{r}{AB}) + \smun{r}{BC}-\tau(\smun{r}{BC}) + \smun{r}{AC}-\tau(\smun{r}{AC}) \] of the direct sum coalgebra structure where we take $r \dsum -\smun{r}{21} \dsum r$ as the quasitriangular structure on the direct sum Lie algebra.  Here, for example, $\smun{r}{AB}=(r^{(1)} \dsum 0 \dsum 0) \tensor (0 \dsum r^{(2)} \dsum 0)$. \end{theorem}

\begin{proof}
Recall the definition of the (Lie algebra) isomorphism $\theta_2$ above, that is, \[ \theta_{2}:T(\Lie{g})  \to \Lie{g} \dsum \Lie{g} \dsum \Lie{g},\ b \dsum g \dsum c \mapsto (b+g) \dsum (b+g-c) \dsum (g-c). \]  We will write $\smun{r}{T}$ for the quasitriangular structure on $T(\Lie{g})$.  From the proposition, we have $\smun{r}{T}= 0\dsum r \dsum 0 + (0 \dsum 0 \dsum f^a) \tensor (e_a \dsum 0 \dsum 0)$ in the direct sum notation, with summation over $a$ understood.  Here $\{ e_a \}$ is a basis of $\bLie{g}$ and $\{ f^a \}$ is a dual basis.  Hence \begin{align*} (\theta_{2} \tensor \theta_{2})(\smun{r}{T}) & = \theta_{2}(0 \dsum r^{(1)} \dsum 0) \tensor \theta_{2}(0 \dsum r^{(2)} \dsum 0) \\ & \quad {} + \theta_{2}(0 \dsum 0 \dsum f^a) \tensor \theta_{2}(e_a \dsum 0 \dsum 0) \\ & = (r^{(1)} \dsum r^{(1)} \dsum r^{(1)}) \tensor (r^{(2)} \dsum r^{(2)} \dsum r^{(2)}) \\ & \quad {} + (0 \dsum -f^a \dsum -f^a) \tensor (e_a \dsum e_a \dsum 0) \end{align*} 

This expression may be simplified as follows.  Label the three copies of $\Lie{g}$ in the direct sum as $\smun{\Lie{g}}{A}$, $\smun{\Lie{g}}{B}$ and $\smun{\Lie{g}}{C}$ and for elements of the tensor product of any two of these, write the appropriate subscripts.  For example, we will write $\smun{a}{AB}$ for $(a_{(1)} \dsum 0 \dsum 0) \tensor (0 \dsum a_{(2)} \dsum 0)$.

We observe that $f^a \tensor e_a$ is precisely the inverse Killing form---or in its quasitriangular form, $2r_{+}$, the symmetric part of $r$.  Hence expanding out the tensor products and rewriting in our subscript notation, we have the following: \begin{align*} (\theta_{2} \tensor \theta_{2})(\smun{r}{T}) & = \smun{r}{AA} + \smun{r}{AB} + \smun{r}{AC} + \smun{r}{BA} + \smun{r}{BB} + \smun{r}{BC} + \smun{r}{CA} + \smun{r}{CB} + \smun{r}{CC} \\ & \quad {} - \smun{(f^a \tensor e_a)}{BA} - \smun{(f^a \tensor e_a)}{BB} - \smun{(f^a \tensor e_a)}{CA} - \smun{(f^a \tensor e_a)}{CB} \end{align*}  However, $f^a \tensor e_a=2r_{+}=r+\tau(r)$ so \begin{eqnarray*} \smun{(f^a \tensor e_a)}{BA} & = & \smun{(r + \tau(r))}{BA} \\ & = & \smun{r^{(1)}}{\!\!B} \tensor \smun{r^{(2)}}{\!\!A} + \smun{r^{(2)}}{\!\!B} \tensor \smun{r^{(1)}}{\!\!A} \\ & = & \smun{r}{BA} + \tau(\smun{r}{AB}) \end{eqnarray*} and similarly for the other terms.  Hence we obtain \begin{align} (\theta_{2} \tensor \theta_{2})(\smun{r}{T}) & = \smun{r}{AA} + \smun{r}{AB} + \smun{r}{AC} + \smun{r}{BA} + \smun{r}{BB} + \smun{r}{BC} + \smun{r}{CA} + \smun{r}{CB} + \smun{r}{CC} \nonumber \\ & \quad {} - (\smun{r}{BA}+\tau(\smun{r}{AB})) - (\smun{r}{BB}+\tau(\smun{r}{BB})) \nonumber \\ & \quad {} - (\smun{r}{CA}+\tau(\smun{r}{AC})) - (\smun{r}{CB}+\tau(\smun{r}{BC})) \nonumber \\ & = \smun{r}{AA} - \tau(\smun{r}{BB}) + \smun{r}{CC} \nonumber \\ & \quad {} + (\smun{r}{AB}-\tau(\smun{r}{AB})) + (\smun{r}{BC}-\tau(\smun{r}{BC})) + (\smun{r}{AC}-\tau(\smun{r}{AC})). \label{rt} \end{align}

Notice first that $\smun{r}{\dsum}=\smun{r}{AA} - \tau(\smun{r}{BB}) + \smun{r}{CC}$ is the direct sum quasitriangular structure on $\Lie{g} \dsum \Lie{g} \dsum \Lie{g}$, choosing the opposite quasitriangular structure for the central copy of $\Lie{g}$.  Now set \[ \chi= \smun{r}{T}-\smun{r}{\dsum}=\smun{r}{AB}-\tau(\smun{r}{AB}) + \smun{r}{BC}-\tau(\smun{r}{BC}) + \smun{r}{AC}-\tau(\smun{r}{AC}).\]  Then $\chi + \smun{\chi}{21}=0$, as is easily seen.  That is, $\chi$ is symmetric and we need only check the identity $(\id \tensor \smun{\delta}{\dsum})\chi + \mathrm{cyclic} + \Lbbracket{\chi}{\chi}= 0$ to see that $\chi$ satisfies the conditions for a cocycle and hence defines a twisting of the direct sum.  This identity follows immediately, however, since we can consider $\chi$ as a sum $\chi = \smun{\chi}{AB} + \smun{\chi}{BC} + \smun{\chi}{AC}$, where $\smun{\chi}{AB}=\smun{r}{AB}-\tau(\smun{r}{AB})$ and similarly for the others.  It is known from the proof that the Drinfel\cprime d double is a twisting of a direct sum (see \cite[Theorem 8.2.5]{FQGT}) that terms of precisely the form $\smun{\chi}{AB}$, etc., satisfy the required identity. \end{proof}

\section{Relationship with the Drinfel\cprime d double}\label{dbreln}

The above description of the triple in Theorem~\ref{twist} is clearly reminiscent of that for the Drinfel\cprime d double (Theorem~\ref{dbtwist}).  More than that, we expect at least one copy of the double to sit inside the triple.  For example, in the bosonisation picture we have $D(\Lie{g})\iso \Lie{g} \ltimesdot \op{\dual{\bLie{g}}}$ and, of course, we defined the triple as $T(\Lie{g})=\dbos{\underline{g}}{g}{\op{\dual{\underline{g}}}}$.  Below, we expand these ideas.

Firstly, let $\Lie{g}$ be a quasitriangular Lie bialgebra, not necessarily factorisable.  We observe that we can write the triple as a matched pair of Lie algebras, as in \cite{MatchedPairs}.  Since we have a Lie algebra structure on the triple, we can break this up as the Lie algebras $\bLie{g}$ and $\Lie{g} \ltimesdot \op{\dual{\bLie{g}}}$ and make a matched pair by actions between them.  Identifying the coalgebra structure on this matched pair, we can rewrite the bialgebra structure of $T(\Lie{g})$ as follows.

\begin{theorem}\label{matpr} Let $\Lie{g}$ be a quasitriangular Lie bialgebra.  Then $T(\Lie{g})$ is a double cross sum Lie algebra and a semidirect Lie coalgebra, written \[ T(\Lie{g})= \bLie{g} \dbxsumrbos (\Lie{g} \ltimesdot \op{\dual{\bLie{g}}}). \] \end{theorem}

\begin{proof} We have a left action of $\Lie{g} \ltimesdot \op{\dual{\bLie{g}}}$ on $\bLie{g}$,
 \begin{align*} & \alpha: (\Lie{g} \ltimesdot \op{\dual{\bLie{g}}}) \tensor \bLie{g} \to \bLie{g} \\ &\alpha((g \dsum c) \tensor b) = \ad{g}{b} - b_{\underline{(1)}}\ip{b_{\underline{(2)}}}{c} \intertext{and a right action of $\bLie{g}$ on $\Lie{g} \ltimesdot \op{\dual{\bLie{g}}}$,}  &\beta:(\Lie{g} \ltimesdot \op{\dual{\bLie{g}}}) \tensor \bLie{g} \to \Lie{g} \ltimesdot \op{\dual{\bLie{g}}} \\ &\beta((g\dsum c) \tensor b)=-c_{\underline{(1)}}\ip{c_{\underline{(2)}}}{b}-2r_{+}^{(1)}\ip{c}{r_{+}^{(2)} \laction b} \end{align*} with $\laction$ being the adjoint action.  Note that these are exactly the terms in the bracket defined on $T(\Lie{g})$ between these two pieces, so we have a matched pair and the Lie bracket on the double cross sum coming from this is exactly that of $T(\Lie{g})$.

We notice that $\bLie{g} \rtimesdot \Lie{g}$ occurs as a sub-Lie bialgebra of $T(\Lie{g})$ and in particular that elements of $\op{\dual{\bLie{g}}}$ do not appear in the cobracket on elements of $\bLie{g}$.  This cobracket is the one obtained by bosonisation, which is by definition a semidirect coalgebra by a Lie coaction.  For $\bLie{g} \rtimesdot \Lie{g}$ the Lie coaction of $\Lie{g}$ on $\bLie{g}$ is $\gamma: \bLie{g} \to \Lie{g} \tensor \bLie{g}$, $\gamma(b)=r^{(2)} \tensor r^{(1)} \laction b$ for $b\in \bLie{g}$.  Then we can extend this coaction to one of $\Lie{g} \ltimesdot \op{\dual{\bLie{g}}}$ on $\bLie{g}$ by letting $\op{\dual{\bLie{g}}}$ coact by zero.   So we can write $T(\Lie{g})=\bLie{g} \rtimesblack (\Lie{g} \ltimesdot \op{\dual{\bLie{g}}})$ as coalgebras. \end{proof}

We can use the isomorphism of $\Lie{g} \ltimesdot \op{\dual{\bLie{g}}}$ with $D(\Lie{g})$ to give a version of this theorem involving the double which is independent of any particular realisation of the double.  The isomorphism is explicitly given by \[ \sigma: D(\Lie{g}) \to \Lie{g} \ltimesdot \op{\dual{\bLie{g}}},\ \sigma(h \dsum d)=(h-r^{(1)}\ip{r^{(2)}}{d}) \dsum d \] for $h\in \Lie{g}$, $d\in \op{\dual{\Lie{g}}}$ and with $D(\Lie{g})=\Lie{g}\dbxsum \op{\dual{\Lie{g}}}$.  Note that the inverse is \[ \inv{\sigma}: \Lie{g} \ltimesdot \op{\dual{\bLie{g}}} \to D(\Lie{g}),\ \inv{\sigma}(g \dsum c)=(g+r^{(1)}\ip{r^{(2)}}{c}) \dsum c \] for $g\in \Lie{g}$, $c\in \op{\dual{\bLie{g}}}$.  That $\sigma$ is a bialgebra isomorphism may be checked along the same lines as the proof in \cite[Example 3.9]{BraidedLie} that $D(\Lie{g}) \iso \dual{\bLie{g}} \rtimesdot \Lie{g}$ (\ie using the opposite conventions).

\begin{corollary}\label{dbxsumrbos} Let $\Lie{g}$ be a quasitriangular Lie bialgebra.  Then $T(\Lie{g})$ is isomorphic to a double cross sum Lie algebra and a semidirect Lie coalgebra \[ T(\Lie{g})\iso \bLie{g} \dbxsumrbos D(\Lie{g}). \] \end{corollary}

\begin{proof} Define the induced actions \[ \hat{\alpha}:D(\Lie{g})\tensor \bLie{g} \to \bLie{g},\ \hat{\alpha}=\alpha \circ (\sigma \tensor \id) \] and \[ \hat{\beta}: D(\Lie{g}) \tensor \bLie{g} \to D(\Lie{g}),\ \hat{\beta}=\inv{\sigma} \circ \beta \circ (\sigma \tensor \id). \]  Explicit expressions may be obtained from Theorem~\ref{matpr}.  These actions give a matched pair $(\bLie{g},D(\Lie{g}))$.  For the coalgebra, we let $\Lie{g}\subset D(\Lie{g})$ coact by $\gamma$ on $\bLie{g}$ as above and let $\op{\dual{\Lie{g}}}$ coact by zero. \end{proof}

We now restrict to $\Lie{g}$ factorisable.  Recall the description above of $T(\Lie{g})$ as a direct sum Lie algebra with twisted coalgebra structure.  Recall also the similar description of $D(\Lie{g})$, as in Theorem~\ref{dbtwist}.  Notice that if we take only the terms in equation~(\ref{rt}) involving the $A$ and $B$ copies, we have $\smun{r}{AB}-\tau(\smun{r}{BB})+\smun{r}{AB}-\tau(\smun{r}{AB})$ which is precisely the quasitriangular structure on the double $D(\Lie{g})$ in the form given by that theorem.  So we can describe the triple in terms of the double as follows.

\begin{corollary}  As Lie bialgebras, $T(\Lie{g}) \iso D(\Lie{g}) \dbxcosum \Lie{g}$ \end{corollary}

\begin{proof}
We see from the proof of Theorem~\ref{twist} and the preceding comments that we have the bialgebra isomorphism $T(\Lie{g})\iso (\Lie{g} \dbxcosum \Lie{g}) \dbxcosum \Lie{g}$. Here, the first double cross cosum is as described in Theorem~\ref{dbtwist} and the twisting in the second is by \[ \smun{\chi}{AC}+\smun{\chi}{BC}= (\smun{r}{AC}-\tau(\smun{r}{AC})) + (\smun{r}{BC}-\tau(\smun{r}{BC})) \] which is a suitable cocycle, as before.  Then the isomorphism of $\Lie{g} \dbxcosum \Lie{g}$ with $D(\Lie{g})$ for $\Lie{g}$ factorisable gives the result.
\end{proof}

\section{Real forms and half-real forms of the triple}\label{realforms}

We now work over $\complex$ and consider real and half-real forms of the triple.  We consider only the factorisable case.  Recall that a real form of a complex Lie algebra $\Lie{g}$ is a choice of basis for $\Lie{g}$ such that all structure constants are real.  There are two particularly natural real forms, namely the split and compact forms, but other non-isomorphic forms too.  See, for example, \cite{FultonHarris} for more on real forms.  The natural basis for a split form gives us a bialgebra over $\reals$ and as the results of the preceding sections hold over any field of characteristic not 2 this case is dealt with. 

We define a half-real form to be a choice of basis with real Lie algebra structure constants and imaginary Lie coalgebra structure constants.  There will, of course, generally be non-isomorphic half-real forms of the same Lie bialgebra.  This concept has been introduced in \cite{MatchedPairs}, as useful in describing Iwasawa decompositions.  Half-real forms are equivalent to real forms but again here we will find them a more useful language.  In particular, the natural choice of basis for compact forms of simple Lie algebras leads us to consider half-real forms.  However, note that a half-real form $(\Lie{u},r)$ is a complex Lie bialgebra, not a bialgebra over $\reals$.  This is because the quasitriangular structure $r$ involves $i$, so $\Lie{u}$ is not quasitriangular over $\reals$.  When $\Lie{g}$ is quasitriangular, we say $r$ is of real type if the symmetric part of $r$, $2r_{+}$, is in fact real.

\begin{lemma} Let $\Lie{g}$ be a complex factorisable Lie bialgebra.  Let $(\Lie{u},r)$ be a half-real form of $\Lie{g}$ with $r$ of real type.  Then the transmutation of $\Lie{u}$, $\bLie{u}$, as described in Definition~\ref{transmut} is a real-real form of the transmutation $\bLie{g}$.  That is, $\bLie{u}$ has real bracket and cobracket structure constants.  Note that $\bLie{u}$ is therefore self-dual.  \end{lemma}

\begin{proof} We are in the case of $r$ real type, \ie $2r_{+}$ real.  Then since by definition $\Lie{u}$ has real bracket structure constants, the braided-Lie cobracket structure constants are real, for we recall that \[ \underline{\delta}x=2r_{+}^{(1)} \tensor \Lbracket{x}{r_{+}^{(2)}} \] for all $x\in \Lie{u}$. Self-duality is ensured by \cite[Example 3.3]{BraidedLie}. \end{proof}

\noindent Therefore we have the following theorem.

\begin{theorem} Let $\Lie{g}$ be a factorisable Lie bialgebra over $\complex$.  Let $(\Lie{u},r)$ be a half-real form of $\Lie{g}$ with real type quasitriangular structure.  Considering $\bLie{u}$ as in the preceding lemma as a complex braided-Lie bialgebra, define the triple of $\Lie{u}$ to be the double-bosonisation $T(\Lie{u})=\dbos{\underline{u}}{u}{\op{\underline{u}}}$.

Then $T(\Lie{u})$ is a half-real form of $T(\Lie{g})$ with real type quasitriangular structure. \end{theorem}

\begin{proof} The brackets on $\Lie{u}$ and $\bLie{u}$ are real by assumption, as are the braided-Lie cobracket $\underline{\delta}$ and the symmetric part $2r_{+}$ of the quasitriangular structure on $\Lie{u}$.  Examining the definitions of the brackets in the triple, we see that these then define a real bracket on $T(\Lie{u})$.

The quasitriangular structure $\smun{r}{T}$ is not real since $r$ on $\Lie{u}$ is not real.  The induced quasitriangular structure on the triple is recalled in Proposition~\ref{doubleqtstr} and may be written as \[ \smun{r}{T}= 0\dsum r \dsum 0 + (0 \dsum 0 \dsum f^a) \tensor (e_a \dsum 0 \dsum 0) \] where $\{ e_a \}$ is a basis of $\bLie{u}$ and $\{ f^a \}$ is a dual basis.  Note that the part $(0 \dsum 0 \dsum f^a) \tensor (e_a \dsum 0 \dsum 0)$ is real.  The dual pairing we are using is the Killing form which is real since $\bLie{u}$ is real as a Lie algebra.  Hence the symmetric part $2(\smun{r}{T})_{+}$ is real: $2r_{+}$ is real and any contribution from $f^{a} \tensor e_{a}$ can only be real.  So $\smun{r}{T}$ is of real type. \end{proof} 

Conversely, if $r$ is not of real type then the bracket on the triple is not real, since $\underline{\delta}$ and $2r_{+}$ are not real.  Then $T(\Lie{u})=\dbos{\underline{u}}{u}{\op{\underline{u}}}$ is not a half-real form of $T(\Lie{g})$.

\section{The triangular case}\label{triang}

Recall that a quasitriangular Lie bialgebra $(\Lie{g},r)$ is said to be triangular if $r$ has zero symmetric part.  Then if $\Lie{b}$ is a $\Lie{g}$-module the associated infinitesimal braiding is also zero: we have \[ \psi(a \tensor b)=2r_{+} \laction (a \tensor b - b\tensor a)=0. \]  So a braided-Lie bialgebra $(\Lie{b},\Lbracket{\ }{\ }_{\Lie{b}},\underline{\delta})$ in the category $\modcat{\Lie{g}}$ is a $\Lie{g}$-covariant bialgebra with $\dform \underline{\delta}=0$.  In particular, we see that this last condition means $\Lie{b}$ is a $\Lie{g}$-covariant unbraided Lie bialgebra in the category.  

We consider the adjoint representation of $\Lie{g}$ in $\modcat{\Lie{g}}$.  The transmutation $\bLie{g}$ as defined in Definition~\ref{transmut} has the adjoint module structure and the Lie bracket of $\Lie{g}$ but the zero braided-Lie cobracket, since $\Lie{g}$ is triangular.  Moreover, this is essentially forced upon us.

\begin{lemma} Let $(\Lie{g},r)$ be a non-Abelian triangular Lie bialgebra and let $\bLie{g}\in \modcat{\Lie{g}}$ be the adjoint representation of $\Lie{g}$, made a Lie algebra in the category by the Lie algebra of $\Lie{g}$.  Let $\underline{\delta}$ be a $\Lie{g}$-covariant cobracket on $\bLie{g}$.  Then $\underline{\delta}=0$. \end{lemma}

\begin{proof} We have $\Lie{g}$-covariance of the coalgebra structure in the form \begin{align} \underline{\delta}(\xi \laction x) & = \xi \laction \underline{\delta}x \intertext{for $\xi\in \Lie{g}$, $x\in \bLie{g}$.  We also have the zero coboundary property for $\underline{\delta}$ as described above.  Explicitly, for $a,b \in \bLie{g}$} \underline{\delta}(\Lbracket{a}{b}) & = \bad{a}{\underline{\delta}b}-\bad{b}{\underline{\delta}a}. \label{cocycle} \end{align}  Here, $\underline{\mathrm{ad}}$ refers to the adjoint action of $\bLie{g}$ on itself, \ie in the category $\modcat{\bLie{g}}$, or its extension to $\bLie{g} \tensor \bLie{g}$.  Define an isomorphism of Lie algebras $\iota : \Lie{g} \to \bLie{g}$ by $\iota = \id$, the identity map.  We can now write $\laction$, the adjoint action of $\Lie{g}$ on $\bLie{g}$, as $\laction: \Lie{g} \tensor \bLie{g} \to \bLie{g}$, $\laction (\xi \tensor x) = \bad{\iota(\xi)}{x}$, or equivalently, $\laction = \underline{\mathrm{ad}} \circ (\iota \tensor \id)$.

So, setting $a=\iota(\xi)$ and $b=x$ in \eqref{cocycle} and writing in terms of $\underline{\mathrm{ad}}$, we have the following two equalities: \begin{align} \underline{\delta}(\bad{\iota(\xi)}{x}) & = \bad{\iota(\xi)}{\underline{\delta}x} \\ \underline{\delta}(\bad{\iota(\xi)}{x}) & = \bad{\iota(\xi)}{\underline{\delta}x}-\bad{x}{\underline{\delta}(\iota(\xi))} \end{align}  Now we see that we have $\bad{x}{\underline{\delta}(\iota(\xi))}=0$ and our choices of $\iota(\xi),x\in \bLie{g}$ were arbitrary.  So we conclude that if $\Lie{g}$ is not Abelian, so that $\underline{\mathrm{ad}}$ and $\iota$ are not identically zero, we must have $\underline{\delta}=0$.  \end{proof}

We now assume $\Lie{g}$ is not Abelian.  Notice that $\bLie{g}$ is now not self-dual, as it has a non-zero bracket and zero braided cobracket.  The dual $\dual{\bLie{g}}$ will have zero bracket and non-zero braided cobracket, namely the (unbraided) Kirillov--Kostant cobracket.  We recall the definition of the triple, using the simplifications we have deduced above.

\begin{definition} Let $\Lie{g}$ be a non-Abelian, finite-dimensional, triangular Lie bialgebra over a field $k$ of characteristic not 2.  Consider the transmutation $\bLie{g}$ as described above.  Define $T(\Lie{g})$ to be the double-bosonisation $\dbos{\underline{g}}{g}{\dual{\underline{g}}}$. \end{definition}

\noindent We have dispensed with the opposite bracket on the dual, as the bracket is zero there.

Examining the bracket from the double-bosonisation, we can write the Lie algebra structure of the triple in this case as follows.

\begin{proposition} Let $\Lie{g}$ be a triangular Lie bialgebra and $T(\Lie{g})$ the triple as defined above.  Then, as a Lie algebra, $T(\Lie{g})$ is isomorphic to a semidirect sum $(\Lie{g} \rtimes_{\mathrm{ad}} \Lie{g}) \ltimes_{\mathrm{coad}} \dual{\Lie{g}}$.  Here $\mathrm{ad}$ and $\mathrm{coad}$ refer to the adjoint and coadjoint actions, respectively, and both parts of $\Lie{g} \rtimes_{\mathrm{ad}} \Lie{g}$ act on $\dual{\Lie{g}}$ by the coadjoint action. \end{proposition}

\begin{proof} This follows immediately from examining the brackets in the double-bosonisation and identifying the non-zero parts.  In particular, we note that $2r_{+}^{(1)}\ip{\phi}{r_{+}^{(2)}\laction x}=-x_{\underline{(1)}}\ip{\phi}{\underline{(2)}}=0$ since $\underline{\delta}x=0$.  We also use $\underline{\delta}\phi=\delta\phi$.  Since $\bLie{g}$ has the same Lie algebra as $\Lie{g}$, on dualisation $\dual{\bLie{g}}$ has the same Lie coalgebra as $\dual{\Lie{g}}$. \end{proof}

\noindent It does not appear that any further simplification of the description of the triple in the triangular case is possible.

\bibliographystyle{halpha}
\bibliography{references}

\vspace{1em}

\noindent\begin{tabular}{ll}
Address:  & 	School of Mathematical Sciences, \\
	  & 	Queen Mary, University of London, \\
	  &	E1 4NS, \\
	  &	United Kingdom. \\
	  &	\\
E-mail:   &	J.Grabowski@qmul.ac.uk \\
Web site: &     http://www.maths.qmul.ac.uk/\webtilde jeg/ \\
	  &	\\
MSC:	  &	17B62
\end{tabular}

\end{document}